\documentclass[reqno,a4paper,final]{amsart}
\usepackage{amssymb,amstext,amsthm,eucal}
\usepackage{bbold}
\usepackage{array}
\usepackage[latin1]{inputenc}
\newcommand{\PO}{\mathrm{PO}} 
\newcommand{\PSU}{\mathrm{PSU}}
\newcommand{\CO}{\mathrm{CO}}   

\newcommand{\OO}{\mathrm{O}}

\newcommand{\SU}{\mathrm{SU}}
\newcommand{\U}{\mathrm{U}}

\newcommand{\RR}{\mathbb{R}}
\newcommand{\CC}{\mathbb{C}}

\newcommand{\eL}{\mathcal{L}}

\newcommand{\eS}{\mathcal{S}}
\newcommand{\eT}{\mathcal{T}}

\newcommand{\eW}{\mathcal{W}}

\newcommand{\eC}{\mathcal{C}}

\newcommand{\eP}{\mathcal{P}}
\newcommand{\eK}{\mathcal{K}}

\newtheorem{THEO}{\bf Theorem}

\newtheorem{COR}{\bf Corollary}
\newtheorem{LM}{\bf Lemma}

\newcommand{\MUNCH}[1]{\relax}
\allowdisplaybreaks[1]
\begin{document}
\begin{sloppypar}
\title{A remark on conformal $\SU(p,q)$-holonomy}
\author{Felipe Leitner}
\address{Institut f{\"u}r Geometrie und Topologie, Universit{\"a}t Stuttgart, Pfaffenwaldring 57,
D-70569 Stuttgart, Germany}
\email{leitner@mathematik.uni-stuttgart.de}
\date{April 2006}

\begin{abstract} If the conformal holonomy group $Hol(\mathcal{T})$ of a simply 
connected space with conformal structure of signature $(2p-1,2q-1)$ is reduced to $\U(p,q)$
then the conformal holonomy is already contained in the special unitary group $\SU(p,q)$. 
We present two different proofs of this statement, one using conformal tractor
calculus and an alternative proof using Sparling's characterisation of Fefferman metrics.  \end{abstract}

\maketitle

\section{Introduction} \label{ab1} 

Any conformal structure on a space $M^n$ of dimension $n\geq 3$ with signature $(r,s)$ 
gives rise to a $|1|$-graded parabolic geometry with canonical Cartan
connection, which is uniquely determined by a normalisation condition on the curvature. This connection 
solves the equivalence problem for conformal geometry.
The standard representation of the orthogonal group $\OO(r+1,s+1)$, 
which double covers the M{\"o}bius group $\PO(r+1,s+1)$, gives rise to a so-called 
standard tractor bundle $\eT(M)$ on $M$ and the
canonical Cartan connection induces a linear connection on this vector bundle.  
This construction allows the definition of the conformal holonomy group 
$Hol(\eT)$ with algebra $\frak{hol}(\eT)$,
which represent important conformal invariants. 

An interesting task in the frame work of conformal geometry is the classification of possible 
conformal holonomy groups and the realisation of such holonomy groups by concrete geometric constructions.
Some progress in this direction was made in the works \cite{Arm05}, \cite{Lei04} and \cite{Lei05a}. 
The probably best known case so far 
is that of (almost) conformally
Einstein spaces. In this situation the conformal holonomy group acts trivially on a certain standard tractor and via a cone
construction a classification of conformal holonomy can be achieved. Obviously, those holonomy groups do not have an irreducible
action on standard tractors. A classical case of an irreducible subgroup of $\OO(2p,2q)$ is the special unitary group $\SU(p,q)$. 
It is well known that conformal structures with holonomy reduction to the subgroup $\SU(p,q)$ are obtained by the classical Fefferman
construction, which assigns to any integrable CR-space of hypersurface type an 
invariant conformal structure on some canonical circle bundle
(cf. \cite{Fef76}, \cite{Lee86}).

We aim to prove here (in two different ways) that if the conformal holonomy group $Hol(\eT)$
of some space $M$ with conformal structure of signature $(2p-1,2q-1)$ is a subgroup
of the unitary group $\U(p,q)$ then the conformal holonomy algebra $\frak{hol}(\eT)$
is automatically contained in 
the special unitary algebra $\frak{su}(p,q)$
(if $M$ is simply connected $Hol(\eT)$ is automatically reduced to $\SU(p,q)$). 
So in the first moment it seems that there is a gap in the list of possible
conformal holonomy groups. However, as one of our proofs shows the reason for this behaviour is directly implicated by the 
normalisation condition for the curvature of the canonical Cartan connection, and therefore, 
the gap should be considered as natural and immediate consequence
in the theory of conformal holonomy. 

We will proceed as follows. In sections \ref{ab2} to \ref{ab4} we will collect and present 
the necessary apparatus to introduce and investigate
conformal holonomy. In particular, we briefly discuss conformal Cartan geometry and tractor 
calculus and make some considerations about conformal
Killing vector fields. In section \ref{ab5} we prove our reduction result 
in terms of tractor calculus (cf. Theorem  \ref{TH1}). In the final
section we give an alternative proof by using Sparlings's characterisation of Fefferman metrics. 
This reasoning allows us to present
a local geometric description of spaces with conformal holonomy reduction to  
$\U(p,q)$ (resp. $\SU(p,q)$) (cf. Corollary \ref{COR100}).

\section{Conformal Cartan geometry} \label{ab2} 

We briefly explain here and in the subsequent section notions and basic facts from conformal Cartan geometry
and tractor calculus. In particular, we will introduce conformal holonomy.
We start with the Cartan geometry. For more detailed material on 
these subjects we refer e.g. to \cite{BEG94}, \cite{CG03}, \cite{CSS97a}, \cite{CSS97b} and \cite{Kob72}.

Let $\RR^{r,s}$ denote the Euclidean space of dimension $n=r+s\geq 3$ and signature $(r,s)$, whereby the 
scalar product is given by the matrix
\[\mathbb{J}=\left(\begin{array}{cc}-I_r& 0\\ 0 & I_{s}\end{array}\right)\ .\]
We denote by $\frak{g}$ the Lie algebra $\frak{so}(r+1,s+1)$ 
of the orthogonal group $\OO(r+1,s+1)$, which acts by the standard representation
on the Euclidean space $\RR^{r+1,s+1}$ of dimension $n+2$ furnished with indefinite scalar product 
\[\langle x,y\rangle=x_-y_++x_+y_-+(x_1,\ldots,x_{n})\ \mathbb{J}\ (y_1,\ldots,y_{n})^\top\ .\]
The Lie algebra $\frak{g}=\frak{so}(r+1,s+1)$ is $|1|$-graded by
\[\frak{g}=\frak{g}_{-1}\oplus\frak{g}_0\oplus \frak{g}_1\ ,\]
where $\frak{g}_0=\frak{co}(r,s)$, $\frak{g}_{-1}=\RR^{n}$ and $\frak{g}_1=\RR^{n*}$.
The $0$-part $\frak{g}_0$ decomposes further to the centre $\RR$ and a semisimple part $\frak{so}(r,s)$, which is the Lie algebra 
of the isometry group of the Euclidean space $\RR^{r,s}$.
We realise the subspaces $\frak{g}_0, \frak{g}_{-1}$ and $\frak{g}_{1}$ of $\frak{g}$ by matrices
of the form
\[\!\!\left(\begin{array}{ccc} 0& 0& 0\\m & 0& 0\\ 0& -m^\top\mathbb{J}& 0
\end{array}\right)\in\frak{g}_{-1}\ ,\quad 
\left(\begin{array}{ccc} -a& 0& 0\\0 & A& 0\\ 0& 0& a
\end{array}\right)\in\frak{g}_0\ ,\quad 
\left(\begin{array}{ccc} 0& l& 0\\0 & 0& -\mathbb{J}\ l^\top\!\!\!\!\\ 0& 0& 0\end{array}\right)\in\frak{g}_{1}\ .\]
The commutators with respect to these matrices are given by
\[ \begin{array}{ll}
{[\ ,\ ]}:\ \frak{g}_0\times\frak{g}_0\to \frak{g}_0\ ,\qquad &
[(A,a),(A',a')]=(AA'-A'A,0)\\[2mm]
{[\ ,\ ]}:\ \frak{g}_0\times \frak{g}_{-1}\to \frak{g}_{-1}\ ,\qquad &
[(A,a),m]=Am+am\\[2mm]
{[\ ,\ ]}:\ \frak{g}_1\times \frak{g}_{0}\to \frak{g}_{1}\ ,\qquad &
[l,(A,a)]=lA+al\\[2mm]
{[\ ,\ ]}:\ \frak{g}_{-1}\times \frak{g}_{1}\to \frak{g}_{0}\ ,\qquad &
[m,l]=(ml-\mathbb{J}\ (ml)^\top\mathbb{J},lm)\ ,\end{array}\qquad \]
where $(A,a)$, $(A',a')\in \frak{so}(r,s)\oplus\RR,\ m\in\RR^{n},\ l\in
\RR^{n*}$.

The space
\[\frak{p}\ :=\ \frak{g}_0\oplus\frak{g}_1\]
is a parabolic subalgebra of $\frak{g}$. We also denote $\frak{p}_+:=\frak{g}_{1}$.
Whilst the grading of $\frak{g}$ is not $\frak{p}$-invariant,
the filtration
\[\frak{g}\ \supset\ \frak{p}\ \supset\ \frak{p}_+\]
is $\frak{p}$-invariant.
The subgroup $P$ of the M{\"o}bius group $G:=\PO(r+1,s+1)$, which consists of those elements whose adjoint action on $\frak{g}$ preserve
the filtration, is a parabolic subgroup with Lie algebra $\frak{p}$. And the subgroup
$G_0$, which preserves the grading of $\frak{g}$, is isomorphic to the group $\CO(r,s)=\OO(r,s)\times\RR_+$ with Lie algebra
$\frak{g}_0=\frak{co}(r,s)$. Moreover, the exponential map restricts to a diffeomorphism from $\frak{p}_+$
onto a normal subgroup $P_+$ of $P$ such that the parabolic $P$ is the semidirect product of 
$G_0$ with $P_+$ (cf. e.g. \cite{CSS97a}).

The map 
\[\begin{array}{ccl}
\partial:\ H\!om(\frak{g}_{-1},\frak{g})&\to& H\!om(\Lambda^2\frak{g}_{-1},\frak{g})\ ,\\[2mm]
\varphi&\mapsto&\partial\varphi\ =\ (\ \  (X,Y)\  \mapsto\ [X,\varphi(Y)]-[Y,\varphi(X)]\ \ )
\end{array}\]
whereby $X,Y\in \frak{g}_{-1}$, 
defines the differential in the cohomology of the abelian Lie algebra $\frak{g}_{-1}$ with values in 
the
representation $\frak{g}$. 
There exists an adjoint map $\partial^*$ to the differential $\partial$ (with respect to a certain inner product on $\frak{g}$),
which is given by
\[\begin{array}{ccl}
\partial^*:\ H\!om(\Lambda^2\frak{g}_{-1},\frak{g})&\to& H\!om(\frak{g}_{-1},\frak{g})\ ,\\[2mm]
\psi&\mapsto& \partial^*\psi\ =\ (\ \ X\in\frak{g}_{-1}\ \mapsto\ \sum_{i=1}^{n}[\eta_i,\psi(\xi_i,X)]\ \ )
\end{array}\]
where $\{\xi_i:\ i=1,\ldots,n\}$ is some basis of $\frak{g}_{-1}$
and $\{\eta_i:\ i=1,\ldots,n\}$ is the corresponding dual basis of $\frak{p}_+$
(with respect to the Killing form of
$\frak{g}$). Alternatively, the adjoint codifferential $\partial^*$ is defined as the negative 
of the dual of the Lie algebra differential for $\frak{p}_+$.  
Both maps $\partial$ and $\partial^*$ are homomorphisms of $\frak{g}_0$-modules (cf. \cite{CSS97b}).

Now let $M^n$ be a $C^\infty$-manifold of dimension $n\geq 3$. Two smooth metrics $g$ and $\tilde{g}$ on $M$
are called conformally equivalent if there exists a $C^\infty$-function $\phi$ such that $\tilde{g}=e^{2\phi}g$.
A conformal structure
$c$ on $M^n$ of signature $(r,s)$, $n=r+s$, is by definition a class of conformally equivalent
$C^\infty$-metrics with signature $(r,s)$. Let $GL(M)$ denote the general linear frame bundle on $M$.  
The choice of a conformal structure $c$ on $M$ is the same thing as a $G_0$-reduction $\mathcal{G}_0(M)$ 
of $GL(M)$. By the process of prolongation we obtain
from $\mathcal{G}_0(M)$ the $P$-principal fibre bundle $\mathcal{P}(M)$, which is a reduction
of the second order linear frame bundle $GL^2(M)$ of $M$ and which again determines the conformal structure $c$ on $M$
uniquely (cf. \cite{CSS97b}, \cite{Kob72}). 

A Cartan connection $\omega$ on a $P$-principal fibre bundle $\mathcal{P}(M)$ is a smooth $1$-form with values
in $\frak{g}$ such that
\begin{enumerate}
\item $\omega(\chi_A)=A\ \ $ for all fundamental fields $\chi_A$, $A\in\frak{p}$,
\item $R_g^*\omega=Ad(g^{-1})\circ\omega\ \ $ for all $g\in P$\ \ and
\item $\omega|_{T_u\mathcal{P}(M)}:T_u\mathcal{P}(M)\to\frak{g}\ \ $ is a linear isomorphism for all
$u\in\mathcal{P}(M)$.
\end{enumerate}
For some arbitrary function $s$ on $\mathcal{P}(M)$ 
we denote by
\[\nabla^\omega_Xs\ :=\ \mathcal{L}_{\omega^{-1}(X)}s\ =\ \omega^{-1}(X)(s)\]
the invariant derivative of $s$ along $X\in\frak{g}_-$ induced by $\omega$.
The curvature $2$-form $\Omega$ of the Cartan connection $\omega$ is given by
\[\Omega\ =\ d\omega+\frac{1}{2}[\omega,\omega]\ .\]
It is $\iota_{\chi_A}\Omega=0$ for all $A\in \frak{p}$, i.e., the curvature $\Omega$ is with respect to insertion of
vertical vectors on $\mathcal{P}(M)$ trivial.
The corresponding curvature function 
\[\kappa:\ \mathcal{P}(M)\ \to\ H\!om(\Lambda^2\frak{g}_-,\frak{g})\]
is defined as
\[\kappa(u)(X,Y):=d\omega(\omega_u^{-1}(X),\omega_u^{-1}(Y))+[X,Y]\]
for $u\in\mathcal{P}(M)$ and $X,Y\in\frak{g}_-$. The curvature function $\kappa$ decomposes 
to $\kappa_{-1}+\kappa_0+\kappa_1$  according to the grading of the target space $\frak{g}$.
The map $\frak{g}_0\to \frak{gl}(\frak{g}_-)$ induced by the adjoint representation
is an inclusion of a subalgebra. Therefore, the $0$-part $\kappa_0$ 
can be seen as a function on $\mathcal{P}(M)$, which takes values in 
$\frak{g}_-^*\otimes\frak{g}_-^*\otimes\frak{g}_-^*\otimes\frak{g}_-$. 
A Cartan connection $\omega$ with the property $\kappa_{-1}\equiv 0$ is called torsion-free (cf. \cite{CSS97a}).

An important fact of conformal geometry is the existence of a canonical 
Cartan connection $\omega_{nor}$ on $\mathcal{P}(M)$, which is 
uniquely determined by the curvature normalisation condition
\[\partial^*\circ\kappa=0\ .\] 
We call $\omega_{nor}$ the normal 
Cartan connection on $\mathcal{P}(M)$ to the space $M$ with conformal structure $c$. 
The condition $\partial^*\circ\kappa=0$ is 
equivalent to 
\[\kappa_{-1}= 0\qquad\mbox{and}\qquad tr\kappa_0:=\sum_{i=1}^n\kappa_0(\xi_i,\cdot)(\cdot)(\eta_i)= 0\ ,\]
i.e., the unique normal Cartan connection $\omega_{nor}$ of conformal geometry is torsion-free and the $0$-part of the 
curvature function is traceless (cf. \cite{CSS97a}, \cite{Kob72}).

The (generalised) Bianchi identity of an arbitrary Cartan connection $\omega$ on $\mathcal{P}(M)$ is expressed by
\[
-\partial\kappa(X,Y,Z)\ =\ \sum_{cycl} (\ \nabla^\omega_Z(\kappa(X,Y))\ +\ \kappa(\kappa_{-1}(X,Y),Z)\ )\ ,
\]
where $\sum_{cycl}$ denotes the sum over all cyclic permutations of the arguments and the map $\partial$
is the Lie algebra differential acting on some $2$-chain
$\psi\in Hom(\Lambda^2\frak{g}_-,\frak{g})$ by
\[ \partial\psi(X,Y,Z)=-\sum_{cycl} [\psi(X,Y),Z] 
\ .\]
For the normal Cartan connection $\omega_{nor}$ on $\mathcal{P}(M)$ the Bianchi identity simplifies to
\begin{eqnarray*}
\sum_{cycl} [\kappa_0(X,Y),Z]&=&0\qquad \mbox{and}\\
\sum_{cycl} \nabla^\omega_Z(\kappa_{0}(X,Y))&=& \sum_{cycl} [\kappa_1(X,Y),Z]\ ,
\end{eqnarray*}
which has much similarity with the first and second Bianchi identity of 
Riemannian geometry (cf. section \ref{ab3} and \cite{CSS97a}).

\section{Some tractor calculus} \label{ab3} 

Let $(M^n,c)$, $n\geq 3$, be a space with conformal structure of signature $(r,s)$ and $\mathcal{P}(M)$ the corresponding $P$-reduction
of $GL^2(M)$ equipped with the normal Cartan connection $\omega_{nor}$. 
The parabolic $P=(\OO(r,s)\times\RR_+)\ltimes \RR^{n*}$ is included in 
$\OO(r+1,s+1)$ by 
\[\begin{array}{cccc}
\iota:& P&\to& \OO(r+1,s+1)\ ,\\[2mm]
&(\beta,\alpha,x)&\mapsto&\left(\begin{array}{ccc}\alpha^{-1}&x&v\\
0&\beta&y\\
0&0&\alpha
\end{array}\right) 
\end{array}\]
whereby $\beta\in \OO(r,s)$, $\alpha>0$ in $\RR$, 
$x\in\RR^{n*}$, $y:=-\alpha\beta\mathbb{J}x^\top$ and $v:=-\frac{\alpha}{2}x\mathbb{J}x^\top$.
This inclusion determines the adjoint action of $P$ on $\frak{g}=\frak{so}(r+1,s+1)$ (such that the filtration
is $P$-invariant). 
Then we define the vector bundle $\mathcal{A}(M)$ over $(M,c)$ by
\[ \mathcal{A}(M):=\mathcal{P}(M)\times_{Ad(P)}\frak{g}\ .\]
We call this associated vector bundle to $\mathcal{P}(M)$ the adjoint tractor bundle on $(M,c)$. Since the adjoint action of $P$ 
preserves the 
filtration of $\frak{g}$, there is a natural filtration on $\mathcal{A}(M)$, which we denote by
\[
\mathcal{A}(M)\ \supset\ \mathcal{A}^0(M)\ \supset\ \mathcal{A}^1(M)\ .
\] 
Thereby, it is $\mathcal{A}^0(M)=\mathcal{P}(M)\times_{Ad(P)}\frak{p}$ 
and the bundle $\mathcal{A}^1(M)=\mathcal{P}(M)\times_{Ad(P)}\frak{p}_+$
is isomorphic to the dual $T^*M$ of the tangent bundle of $M$. 
In fact, the dual tangent bundle $T^*M$ is canonically contained in the 
adjoint tractor bundle $\mathcal{A}(M)$. Moreover, the quotient $\mathcal{A}(M)/\mathcal{A}^0(M)$ is 
naturally identified with the tangent bundle $TM$. For a section $A\in\Gamma(\mathcal{A}(M))$ we denote by $\Pi(A)\in\frak{X}(M)$ 
the corresponding projection to vector fields on $M$. 

By restriction of the standard representation of $\OO(r+1,s+1)$ to $P$ via the inclusion $\iota$
we obtain a representation of $P$ on the Euclidean space $\RR^{r+1,s+1}$, 
which also preserves the scalar product $\langle\cdot,\cdot\rangle$ (cf. 
section \ref{ab2}).
We use this $P$-representation to define the standard tractor bundle $\mathcal{T}(M)$ over the space $(M,c)$ by 
\[\mathcal{T}(M):=\mathcal{P}(M)\times_{P}\RR^{r+1,s+1}\ .\]
The standard tractor bundle is equipped with an invariant scalar product 
$\langle\cdot,\cdot\rangle$ and it admits a filtration as well, which 
we denote by
\[
\mathcal{T}(M)\ \supset\ \mathcal{T}^0(M)\ \supset\ \mathcal{T}^1(M)\ .
\]
Thereby, it is  $\mathcal{T}^1(M)=\mathcal{P}(M)\times_{P}\RR$, a trivial real line bundle, which is isomorphic to the density bundle
$\mathcal{E}[-1]$ induced from the representation of conformal weight $1/2$.
(If $\mathcal{R}\subset S^2(T^*M)$ denotes the ray subbundle of metrics $g$ in $c$ with
$\RR_+$-action $\mathcal{R}\times\RR_+\to\mathcal{R}$ given by $(g,s)\mapsto s^2\cdot g$ then the densities $\mathcal{E}[w]$
are defined  as associated line bundles by $\mathcal{E}[w]:=\mathcal{R}\times_{\rho^w}\RR$, where $\rho^w:a\in\RR_+\mapsto a^{-w/2}
\in Aut(\RR)$ (cf. e.g. \cite{CG03}).) 
Moreover, since the standard $\frak{g}$-action on $\RR^{r+1,s+1}$
is compatible with the $P$-action, there exists a natural equivariant action of $\mathcal{A}(M)$ on $\mathcal{T}(M)$, 
which we denote by
\[\bullet:\mathcal{A}(M)\otimes\mathcal{T}(M)\to\mathcal{T}(M)\ .\]
In this respect we can understand $\mathcal{A}(M)$ as a subbundle of the endomorphism bundle $End(\mathcal{T}(M))$.

The tractor bundles $\mathcal{T}(M)$ and $\mathcal{A}(M)$ over $(M,c)$ 
enjoy the existence of naturally  defined covariant derivatives.
To introduce these covariant derivative we extend the normal Cartan connection $\omega_{nor}$ to a principal fibre bundle
connection in the standard way as follows. We set 
\[\tilde{\mathcal{G}}(M):=\mathcal{P}(M)\times_{\iota}\OO(r+1,s+1)\ .\]
This $\OO(r+1,s+1)$-bundle together with the natural right action of $\OO(r+1,s+1)$ is a principal fibre bundle over $(M,c)$.
The normal Cartan connection $\omega_{nor}$ on $\mathcal{P}(M)$ extends by right translation with the action of $\OO(r+1,s+1)$ 
to a principal fibre bundle connection on $\tilde{\mathcal{G}}(M)$, which we denote by $\omega_{nor}$ as well.
The tractor bundles $\mathcal{T}(M)$ and $\mathcal{A}(M)$ are associated vector bundles 
to the principal fibre bundle $\tilde{\mathcal{G}}(M)$. In fact, it holds
\[
\mathcal{T}(M)=\tilde{\mathcal{G}}(M)\times_{\OO(r+1,s+1)}\RR^{r+1,s+1}\qquad
\mbox{and}\qquad\mathcal{A}(M)=\tilde{\mathcal{G}}(M)\times_{Ad}\frak{g}
\]   
and the principal fibre bundle connection $\omega_{nor}$ induces in the usual manner linear connections on these vector bundles. 
We denote these connections by 
\[
\nabla^{nor}:\Gamma(\mathcal{T}(M))\to \Gamma(T^*M\otimes\mathcal{T}(M))\quad\mbox{and}\quad
\nabla^{nor}:\Gamma(\mathcal{A}(M))\to \Gamma(T^*M\otimes\mathcal{A}(M)).
\]

The curvature $\Omega$ of the normal Cartan connection $\omega_{nor}$ on
$\mathcal{P}(M)$ is by definition $P$-equivariant and  the equivariant extension of $\Omega$ 
to $\tilde{\mathcal{G}}(M)$ equals the curvature of the principal fibre bundle connection $\omega_{nor}$.
In particular,
$\Omega$ can be understood as a smooth section in the adjoint tractor bundle $\mathcal{A}(M)$ and 
it acts via $\bullet$ as 
the curvature operator of $\nabla^{nor}$, i.e., it holds
\[
\nabla^{nor}_X\nabla^{nor}_Yt-\nabla^{nor}_Y\nabla^{nor}_Xt-\nabla^{nor}_{[X,Y]}t=\Omega(X,Y)\bullet t\]
for all standard tractors $t\in\Gamma(\mathcal{T}(M))$ and $X,Y\in\frak{X}(M)$.

The covariant derivative $\nabla^{nor}$ induces a parallel displacement on the standard tractor bundle $\mathcal{T}(M)$ over $(M,c)$.
In particular, the parallel displacement along loops in $(M,c)$ generates a Lie subgroup of the structure group $\OO(r+1,s+1)$.
We denote this Lie group by $Hol(\mathcal{T})$ and call it the conformal holonomy group of $(M,c)$.
Moreover, we denote by $\frak{hol}(\mathcal{T})$ the holonomy algebra to $Hol(\mathcal{T})$. 
Both objects $Hol(\mathcal{T})$ and $\frak{hol}(\mathcal{T})$
are conformal invariants naturally attached to the underlying space $(M,c)$ (cf. \cite{Arm05}, \cite{Lei04}, \cite{Lei05a}). 

Now we choose some arbitrary metric $g$ in the conformal class $c$ on $M$
and describe the tractor bundles $\mathcal{T}(M)$ and $\mathcal{A}(M)$ and the curvature $\Omega$ of the 
normal Cartan connection $\omega_{nor}$ with respect to this metric $g$. First of all, we note that $g$ corresponds in a unique way
to a $G_0$-equivariant section \[\rho_g:\ \mathcal{G}_0(M)\ \to\ \mathcal{P}(M)\ ,\] i.e., $\rho_g$ is a $G_0$-equivariant lift 
in the prolongation by $P_+$ from the $G_0$-reduction
$\mathcal{G}_0(M)$ of $GL(M)$ to the $P$-reduction $\mathcal{P}(M)$ of the second order frame bundle $GL^2(M)$. 
With the help of the lift $\rho_g$ we can reduce the structure group of the tractor bundles $\mathcal{T}(M)$ and $\mathcal{A}(M)$
to $G_0=\CO(r,s)$ and we obtain the following identifications of vector bundles,
\[
\mathcal{T}(M)\ \cong_g\ \mathcal{G}_0(M)\times_{\iota(G_0)}\RR^{r+1,s+1}\qquad\mbox{and}\qquad\mathcal{A}(M)\ \cong_g\ 
\mathcal{G}_0(M)\times_{Ad\circ\iota(G_0)}\frak{g}\ .
\]
It is important to note that these identifications (denoted by $\cong_g$) depend strongly on 
the choice of the metric $g$ (resp. the corresponding lift $\rho_g$). 

Moreover, the restricted representations of $G_0$ on $\RR^{r+1,s+1}$ and $\frak{g}$ are not any longer indecomposable. 
In fact, $\RR^{r+1,s+1}$
decomposes as $G_0$-module into the direct sum $\RR(1)\oplus \RR^{r,s}\oplus \RR(-1)$, whereby $\RR(w)$ denotes the 
real line of conformal weight $-w/2$, and $\frak{g}$ decomposes into the 
grading $\frak{g}_{-1}\oplus\frak{g}_0\oplus\frak{g}_1$. Accordingly, via 
the lift $\rho_g$ we obtain identifications of $\mathcal{T}(M)$ and $\mathcal{A}(M)$ with graded vector bundles:   
\[\mathcal{T}(M)\ \cong_g\ \mathcal{E}[1]\oplus TM[-1]\oplus\mathcal{E}[-1]\!\qquad\mbox{and}\!\qquad
\mathcal{A}(M)\ \cong_g\ TM\oplus \frak{co}(TM)\oplus T^*M,\]
where $TM[-1]:= TM\otimes\mathcal{E}[-1]$ has zero conformal weight.
With respect to the scalar product $\langle\cdot,\cdot\rangle$ on $\mathcal{T}(M)$ 
the subbundles $\mathcal{E}[1]$ and $\mathcal{E}[-1]$
are lightlike, i.e., it holds $\langle u,u\rangle=0$ for any $u\in \mathcal{E}[1]$ resp. $u\in \mathcal{E}[-1]$.
The restriction of $\langle\cdot,\cdot\rangle$ to $TM$ in $\mathcal{T}(M)$ gives rise to the metric $g$ again. 

With a choice of a metric $g$ and the corresponding gradings on $\mathcal{T}(M)$ and $\mathcal{A}(M)$ 
we can express tractors $t\in\Gamma(\mathcal{T}(M))$ and $A\in\Gamma(\mathcal{A}(M))$ 
as triples $t=(b,\tau,d)$ and $A=(\xi,\varphi,\eta)$.
In particular, the curvature $\Omega\in\Gamma(\mathcal{A}(M))$ of $\nabla^{nor}$ 
on $\mathcal{T}(M)$ decomposes to $(\Omega_{-1},\Omega_0,\Omega_{1})$.
The normalisation condition for $\Omega$ implies that $\Omega_{-1}$ and the trace 
of $\Omega_{0}\in\Omega^2(M,\frak{co}(TM))$ vanish identically.
We note that the three components $\Omega_{-1}$, $\Omega_0$ and  $\Omega_1$ have well known 
interpretations in terms of the metric $g$ and its curvature expressions.
In fact, the $2$-form $\Omega_{-1}$ with values in $TM$ is the torsion 
of the Levi-Civita connection $\nabla^g$ to $g$, which is known to be zero. 
The $2$-form $\Omega_{0}$ with values in the skew-adjoint endomorphism $\frak{co}(TM)$ is known to be the trace-free part of 
the Riemannian curvature tensor of $g$, which 
is called the Weyl tensor and is usually denoted by $\eW$. Eventually, the $1$-part $\Omega_{1}$ 
is equal to the Cotton-York tensor $\eC$ of $g$, which is given by
\[\eC(X,Y)\ :=\ (\nabla^g_X\eK)(Y)-(\nabla^g_Y\eK)(X)\ ,\]
whereby $Ric^g$ denotes the Ricci curvature of $g$, $scal^g$ the scalar curvature and  
\[\eK\ :=\ \frac{1}{n-2}\left(\frac{scal^g}{2(n-1)}g-Ric^g\right)\]
is the Schouten tensor.
The Weyl tensor $\eW$ of $g$ satisfies the first Bianchi identity 
\[
\sum_{cycl}\eW(X,Y)Z\ =\ 0\ ,
\]
which is a direct consequence of the generalised Bianchi identity for the curvature function 
$\kappa$ of the normal Cartan connection $\omega_{nor}$
(cf. section \ref{ab2}). As a $(4,0)$-tensor the Weyl tensor $\eW$ has also the symmetry property
$\eW(X,Y,Z,U)=\eW(Z,U,X,Y)$ for $X,Y,Z,U\in TM$ and the cyclic sum over the last three arguments of $\eW(X,Y,Z,U)$ is zero as well.
Since the torsion part of $\Omega$ vanishes, the Weyl tensor $\eW$ 
expresses a conformal invariant of $(M,c)$, whereas 
the Cotton-York tensor $\eC$ itself is not conformally invariant.

Finally, in this section, we introduce some conventions about calculations with tractors.
First, as we explained above it is $\eT(M)\cong_g\mathcal{E}[1]\oplus TM[-1]\oplus \mathcal{E}[-1]$ 
with respect to the choice of a metric
$g$. A standard tractor $t\in\mathcal{T}(M)$ is then given by a triple $(b,\tau,d)$. 
However, the metric $g$ does not only induce a reduction of $\eP(M)$ to the group 
$G_0=\CO(r,s)$, but even gives rise
to a $\OO(r,s)$-reduction $\mathcal{G}_{ss}(M)$ of $GL(M)$. 
Via this reduction we obtain unique identifications for all density bundles 
$\mathcal{E}[w]$ with the weightless trivial real line bundle over $M$, which we simply denote by $\RR$. In particular,
we have a unique identification
\[\eT(M)\ \cong_g\ \RR\ \oplus\ TM\ \oplus\ \RR\]
and a standard tractor $t\in\mathcal{T}(M)$ is given by
a triple $(b,\tau,d)$, where $b,d\in\RR$ and $\tau$ is simply a tangent vector. We will often 
express this triple as a column vector 
\[ t=\left(\begin{array}{c}d\\ \tau\\ b\end{array}\right)\ .\]

Accordingly, we will write an adjoint tractor $A=(\xi,\varphi,\eta)\in\mathcal{A}(M)$ with respect to the grading induced by $g$ 
as a matrix of the form
\[
A\ =\ \left(\begin{array}{ccc}
-\varphi_c&\eta&0\\
\xi&\varphi_{ss}&-\eta^\flat\\
0&-g(\xi,\cdot)&\varphi_c
\end{array}\right)\ ,
\]
whereby $\varphi=\varphi_{ss}+\varphi_c\cdot id|_{TM}$, $\varphi_{ss}$ is skew-symmetric 
and $\eta^\flat$ denotes the dual tangent vector to $\eta$ with respect to $g$. 
In particular, the curvature $\Omega$ of $\nabla^{nor}$
is presented
with respect to a metric $g$ by the matrix
\[
\left(\begin{array}{ccc}
0&\Omega_{1}&0\\
\Omega_{-1}&\Omega_0&-\Omega_{1}^\flat\\
0&-g(\Omega_{-1},\cdot)&0
\end{array}\right)\ =\ 
\left(\begin{array}{ccc}
0&-\eC&0\\
0&\eW&\eC^\flat\\
0&0&0
\end{array}\right)\ .
\]
The action $\bullet$ of an adjoint tractor $A$ on a standard tractor $t$ can then be expressed by matrix multiplication:
\[
A\bullet t\ =\ \left(\begin{array}{ccc}
-\varphi_c&\eta&0\\
\xi&\varphi_{ss}&-\eta^\flat\\
0&-g(\xi,\cdot)&\varphi_c
\end{array}\right)\left(\begin{array}{c}d\\ \tau\\ b\end{array}\right)\ .
\]
  
\section{Conformal Killing vector fields as adjoint tractors} \label{ab4} 

In this section we recall a relation of conformal Killing vector fields and adjoint tractors, which satisfy
a certain tractor equation. In particular, we will establish Lemma \ref{LM0}, which will be crucial in the following sections 
and which was proved (in a much more general setting) in \cite{Cap05} (cf. also \cite{Lei05a}).    

Let $(M,c)$ be a space with conformal structure. A vector field $V\in\frak{X}(M)$ is called a conformal Killing vector 
field if for some metric $g\in c$ and some function $\lambda$ it holds \[\eL_Vg=\lambda\cdot g\ ,\] 
where $\eL$ denotes the Lie derivative.
In fact, if this condition is satisfied by $V$ then it exists for any metric $\tilde{g}\in c$ some function $\tilde{\lambda}$
such that $\eL_V\tilde{g}=\tilde{\lambda}\cdot \tilde{g}$. 
Let us consider now the $P$-reduction $\eP(M)$ of the second order  
frame bundle $GL^2(M)$ with projection $\pi$ to $(M,c)$. The $2$-jet $j_2(V)$ of a conformal Killing vector field $V$ 
can be interpreted in a natural way as a vector field 
on $\eP(M)$, which then satisfies the condition \[\mathcal{L}_{j_2(V)}\omega_{nor}=0\ .\] 
On the other hand, any vector field $Q\in\frak{X}(\eP(M))$ 
with the property $\mathcal{L}_{Q}\omega_{nor}=0$ projects uniquely to a conformal 
Killing vector field $\pi_*Q$ on $(M,c)$ (and the $2$-jet of $\pi_*Q$
equals $Q$ again). 

Furthermore, let $Q=j_2(V)$ be the $2$-jet of a conformal Killing vector field $V$. 
The evaluation of $Q$ with the normal Cartan connection $\omega_{nor}$
produces the function
\[\begin{array}{cccl}A_Q:& \eP(M)&\to&\ \frak{g}\ .\\[1mm]
& u&\mapsto& \omega_{nor}(Q(u))
\end{array}
\]
This function $A_Q$ is $P$-equivariant on $\mathcal{P}(M)$. Hence $A_Q$ can be interpreted as a section in the adjoint tractor bundle
$\mathcal{A}(M)$. The quotient $\Pi(A_Q)\in\frak{X}(M)$ of this adjoint tractor equals $V=\pi_*Q$. 

\begin{LM} \label{LM-1} (cf. \cite{Cap05})
Let $Q\in\frak{X}(\eP(M))^P$ be a $P$-invariant vector field and let $A_Q$ be 
the corresponding section in the adjoint tractor bundle $\mathcal{A}(M)$ (via $\omega_{nor}$).
\begin{enumerate}
\item\label{z1}
The condition $\mathcal{L}_{Q}\omega_{nor}=0$ for the vector $Q$ is equivalent 
to the condition \[\nabla^{nor}_XA_Q=-\Omega(\Pi(A_Q),X)\qquad\mbox{for\ all}\
\ X\in TM\] 
on the adjoint tractor $A_Q\in\Gamma(\mathcal{A}(M))$. 
\item
If (one of) the conditions of (\ref{z1}) are satisfied then $\pi_*Q=\Pi(A_Q)\in\frak{X}(M)$ is a conformal Killing vector field.  
\end{enumerate}
\end{LM}

As next we consider the tractor equation \[\nabla^{nor}A=0\] for a section $A\in\Gamma(\mathcal{A}(M))$.
The following result is not simply a matter of course.
\begin{LM} \label{LM0} (cf. \cite{Cap05}) Let $A\in\Gamma(\mathcal{A}(M))$ be a $\nabla^{nor}$-parallel adjoint tractor. 
Then it holds \[
\Omega(\Pi(A),\cdot)=0\ .
\]
In particular, the quotient $\Pi(A)$ 
is a conformal Killing vector field (and $A$ corresponds via $\omega_{nor}$ to the $2$-jet of $\Pi(A)$). 
\end{LM}
We note that Lemma \ref{LM0} is true for all regular normal parabolic geometries. In particular, 
since $\frak{g}=\frak{so}(r+1,s+1)$ is $|1|$-graded,
it holds for normal conformal geometries
$(\mathcal{P}(M),\omega_{nor})$.

Finally, we present with respect to any metric $g\in c$ the second order differential operator, which assigns to some conformal
Killing vector field $V$ on $M$ the adjoint tractor $A_Q$, which corresponds to the $2$-jet $Q:=j_2(V)$. It holds
\[
A_Q\ \simeq_g\ (V,\nabla^gV,\mathcal{D}^g(V))\ =\ \left(\begin{array}{ccc}
-\frac{2}{n}div^gV&\mathcal{D}^g(V)&0\\[1.5mm]
V&asym \nabla^gV&-\mathcal{D}^g(V)^\flat\\[1.5mm]
0&-g(V,\cdot)&\frac{2}{n}div^gV
\end{array}\right)\ ,
\]
whereby $div^gV:=tr_g(\nabla^g_\cdot V,\cdot)$ is the divergence, 
the covariant derivative $\nabla^gV$ splits into the anti-symmetric part $asym\nabla^gV$ 
and the symmetric part $div^gV\cdot id|_{TM}$ and
\[
\mathcal{D}^g\ :=\ \frac{1}{n-2}\left(\ tr_g\nabla^2\ +\ \frac{scal^g}{2(n-1)}\ \right)
\]
is the Bochner-Laplacian with a curvature normalisation (cf. \cite{Lei05b}).
This differential operator acting on vector fields and mapping to adjoint tractors is well known in 
the theory of the construction of BGG-sequences as splitting operator. 
Splitting operators exist in much more general fashion (cf. \cite{CSS01}).

\section{Conformal holonomy reduction} \label{ab5}

In section \ref{ab3} we introduced the conformal holonomy group $Hol(\mathcal{T})$ and the holonomy algebra $\frak{hol}(\mathcal{T})$,
which describe conformal invariants for any space $M$ with conformal structure $c$. In this section, we will assume that there exists 
a section $J$ of the adjoint tractor bundle $\mathcal{A}(M)$ on $(M^n,c)$ of dimension $n=2m\geq 4$ with signature 
$(2p-1,2q-1)$, which acts as complex structure on the standard tractors 
$\mathcal{T}(M)$ and which is parallel with respect to the covariant derivative $\nabla^{nor}$, i.e., it holds 
\[J^2:=J\bullet J=-id|_{\mathcal{T}(M)}\in\Gamma(End(\mathcal{T}(M))\qquad\mbox{and}\qquad 
\nabla^{nor}J=0\ .\] The existence of such a complex structure $J$ on 
$\mathcal{T}(M)$ is equivalent (by the very definition of holonomy) 
to the fact that the holonomy group
$Hol(\mathcal{T})$ of $(M,c)$ is contained in the unitary group $\U(p,q)$.
However, we aim to show here in terms of tractor calculus that the 
existence of such $J$ 
implies already that the holonomy algebra $\frak{hol}(\mathcal{T})$ is even reduced to the special unitary group $\frak{su}(p,q)$! 
The reason for this reduction is essentially the normalisation 
condition $\partial^*\circ\kappa=0$ on the curvature of the 
canonical Cartan connection $\omega_{nor}$. 
We will explain later (cf. section \ref{ab6}) that our result follows alternatively
from a result of
C.R. Graham about Sparling's characterisation of Fefferman metrics (cf. \cite{Spa85}, \cite{Gra87}). 

We start our reasoning with an observation about complex structures in $\frak{so}(2p,2q)$. 
\begin{LM} (cf. \cite{Lei05b}) \label{LM1} Let 
\[\beta\ =\ 
\left(\begin{array}{ccc} -a& l& 0\\m & A& -\mathbb{J}l^\top\\ 0& -m^\top\mathbb{J}& a
\end{array}\right)
\]
be a matrix in $\frak{g}=\frak{so}(2p,2q)$. Then the property $\beta^2=-id$ is equivalent to the following conditions on $m,
-\mathbb{J}l^\top\in\frak{g}_{-1}$ and $A\in\frak{so}(\mathbb{J})$:
\begin{enumerate}
\item $m$ and $-\mathbb{J}l^\top$ are lightlike eigenvectors of $A$ to the eigenvalue $a$,
\item the scalar product of $m$ with $-\mathbb{J}l^\top$ equals $1+a^2$ and 
\item $A^2$ restricted to the complement $span\{m,-\mathbb{J}l^\top\}^\bot$ in $\frak{g}_{-1}$ is equal to $-id$.  
\end{enumerate}
\end{LM}

{\bf Proof.} It is
\[
\beta^2=\left(\
\begin{array}{ccc}
a^2+lm& -al+lA&-l\mathbb{J}\ {}^tl\\[2mm]
-am+Am&ml+A^2+\mathbb{J}\ {}^tl{}^tm\mathbb{J}& -A\mathbb{J}\ {}^tl-a\mathbb{J}\ {}^tl\\[2mm]
-{}^tm\mathbb{J}m&-{}^tm\mathbb{J}A-a{}^tm\mathbb{J}& {}^tm{}^tl+a^2
\end{array}\right)\quad .
\]
From this result about the matrix square of $\beta$ the statement of Lemma \ref{LM1} becomes obvious.
\hfill$\Box$\\

Let $J\in\mathcal{A}_x(M)$ be a complex structure of $\mathcal{T}_x(M)$ at some point $x\in M$, 
i.e., $J^2=J\bullet J=-id|_{\mathcal{T}_x(M)}$,
and let $g \in c$ be an arbitrary metric on $M$. With respect to the grading of $\mathcal{A}_x(M)$ induced by $g$ we can write
the complex structure $J$ according to Lemma \ref{LM1} as a matrix
\[
J\ =\
\left(\begin{array}{ccc} -J_c& J_\eta& 0\\ j & J_{ss}& -J_\eta^\flat\\ 0& -g(j,\cdot)& J_c
\end{array}\right)\ ,
\]
whereby $j, -J_\eta^\flat\in T_xM$ are lightlike $J_c$-eigenvectors of $J_{ss}$ with $-g(j,J_\eta^\flat)=1+J_c^2$.
We define the subspace $W(J,g)$ of $T_xM$ as the $g$-orthogonal complement to $span\{j,J_\eta^\flat\}$ in $T_xM$. Then it holds
\[\langle\cdot,\cdot\rangle|_{W(J,g)}\ =\ g|_{W(J,g)} \qquad\mbox{and}\qquad J|_{W(J,g)}\ =\ J_{ss}|_{W(J,g)}\ ,\]
whereby the restriction of the endomorphism $J_{ss}$ acts as $g$-orthogonal complex structure on $W(J,g)$. 
We note that if $u_-$ and $u_+$ are generating elements of 
$\mathcal{E}_x[1]\cong_g\RR$ resp. $\mathcal{E}_x[-1]\cong_g\RR$ with $\langle u_-,u_+\rangle=1$
then we have
\[
W(J,g)\ :=\ span\{j,J_\eta^\flat\}^{\bot_g}\ 
=\ span\{u_-,Ju_-,u_+,Ju_+\}^{\bot_{(\mathcal{T},\langle\cdot,\cdot\rangle)}}\ \subset\ T_xM\ \subset \mathcal{T}_x(M).
\]
In this situation we can choose a complex basis $\{e_\alpha:\ \alpha=1,\ldots, m-1\}$ 
of $W(J,g)$ such that $\{e_\alpha,Je_\alpha:\ \alpha=1,\ldots, m-1\}$
is an orthogonal basis of $W(J,g)$ and  
\[
\{\ u_-\ ,\ Ju_-\ ,\ u_+\ ,\ Ju_+\ ,\ e_1\ ,\ Je_1\ ,\ \ldots\ ,\ e_{m-1}\ ,\ Je_{m-1}\ \}
\] 
is a basis of the space $\mathcal{T}_x(M)$ of standard tractors at $x\in M$. We call the complex basis of the form 
\[\mathcal{B}(J,g)\ :=\ \{\ u_-\ ,\ u_+\ ,\ e_1\ ,\ \ldots\ ,\ e_{m-1}\ \}\]
a $(J,g)$-adapted basis of $\mathcal{T}_x(M)$. 

Now let $(M^n,c)$ be a space of dimension $n=2m\geq 4$ with conformal structure of signature $(2p-1,2q-1)$ 
and a $\nabla^{nor}$-parallel complex structure $J\in\Gamma(\mathcal{A}(M))$.
The holonomy group $Hol(\mathcal{T})$ is reduced to $\U(p,q)$.  
Let us denote by \[\mathcal{T}^\CC(M)=\mathcal{T}(M)\otimes\CC\] the complexified standard tractor bundle. 
We extend the complex structure $J$ on $\mathcal{T}(M)$ to a $\CC$-linear 
complex structure on the complexification $\mathcal{T}^\CC(M)$,
which we denote again by $J$. The bundle $\mathcal{T}^\CC(M)$ decomposes into the direct 
sum $\mathcal{T}_{10}\oplus\mathcal{T}_{01}$, whereby
$\mathcal{T}_{10}$ denotes the $i$-eigenspace of $J$ and $\mathcal{T}_{01}$ is the complex conjugate. The determinant bundle 
\[\mathcal{S}:=\Lambda^{m+1*}(\mathcal{T}_{10})\]
is a complex line bundle on $M$. We call $\mathcal{S}$ the canonical complex 
line tractor bundle of $(\mathcal{T}(M),J)$. 
(If we denote by $\tilde{\mathcal{U}}(M)$ the $\U(p,q)$-reduction 
of the $\OO(2p,2q)$-principal fibre bundle $\tilde{\mathcal{G}}(M)$
induced by $J$ then the canonical complex line tractor bundle is given 
by $\mathcal{S}=\tilde{\mathcal{U}}(M)\times_{det_{\CC}^{-1}}\CC$.)     
  
The principal fibre bundle connection $\omega_{nor}$ induces on $\mathcal{S}$ a covariant derivative, 
which we denote again by $\nabla^{nor}$.
We also denote by $\Omega_{\mathcal{S}}$ the (conformal) curvature of this connection.
\begin{LM}\label{LM2}
Let $\mathcal{S}$ be the canonical complex line tractor bundle of $\mathcal{T}(M)$ 
with $\nabla^{nor}$-parallel complex structure $J$.
Then the curvature $\Omega_\mathcal{S}$ on $\mathcal{S}$ induced by the normal Cartan connection 
$\omega_{nor}$ vanishes identically.
\end{LM} 

{\bf Proof.}
We aim to compute the curvature $\Omega_{\mathcal{S}}$ on $\mathcal{S}$. First of all, 
we remark that with Lemma \ref{LM0} and the assumptions we know that 
\[ \Omega(\Pi(J),\ \cdot\ )\ =\ 0\ .\]
We set $j:=\Pi(J)$ and with respect to an arbitrary metric $g\in c$ we can conclude that 
$\iota_{j}\eW=0$ and $\iota_{j}\eC=0$.

Now let $\{E_\alpha:\ \alpha=1,\ldots,m+1\}$ be a local complex frame of 
$(\mathcal{T}(M),J)$
such that $\{E_\alpha,JE_\alpha:\ \alpha=1,\ldots,m+1\}$ is a local orthonormal frame of 
$(\mathcal{T}(M),\langle\cdot,\cdot\rangle)$. 
Then we denote \[\varrho_\alpha\ :=\ \frac{1}{\sqrt{2}}\langle\ \cdot\ ,\ E_\alpha+iJE_\alpha\ \rangle\]  
and the $(m+1)$-form \[\varrho\ :=\ 
\varrho_1\wedge\ldots\wedge\varrho_{m+1}\]
is a local complex tractor volume form on $(M,c)$, i.e., a local section in $\mathcal{S}$. It holds
\[\Omega_{\mathcal{S}}(X,Y)\circ\varrho\ =\ -i\sum_{\alpha=1}^{m+1} 
\langle E_\alpha, E_\alpha\rangle\cdot \langle \Omega(X,Y)\bullet E_\alpha, JE_\alpha\rangle\cdot 
\varrho
\]
for all $X,Y\in TM$. In particular, this expression proves that \[\Omega_\mathcal{S}(j,\ \cdot\ )\ =\ 0\ . \]

As next we reformulate above expression for $\Omega_{\mathcal{S}}$ 
with respect to an arbitrary metric $g\in c$ on $M$ and a local $(J,g)$-adapted 
frame $\mathcal{B}(J,g)=\{u_-,u_+,e_1,\ldots,e_{m-1}\}$ with $\varepsilon_\alpha:=g(e_\alpha,e_\alpha)$. 
It holds
\[
\Omega_\mathcal{S}(X,Y)\ =\ -2i\cdot\langle \Omega(X,Y)\bullet u_-,Ju_+\rangle\ 
-\ i\sum_{\alpha=1}^{m-1}\varepsilon_\alpha\cdot \langle \Omega(X,Y)\bullet 
e_\alpha, Je_\alpha\rangle\ , 
\]
which is a purely imaginary number for $X,Y\in TM$. We remember that the curvature $\Omega$ 
has with respect to $g$ the matrix form 
\[\left(\begin{array}{ccc}
0&\Omega_1&0\\
\Omega_{-1}&\Omega_0&-\Omega_1^\flat\\
0&-g(\Omega_{-1},\cdot)&0
\end{array}\right)
\ =\
\left(\begin{array}{ccc}
0&-\eC&0\\
0&\eW&\eC^\flat\\
0&0&0
\end{array}\right)\ .
\]
We obtain 
\begin{eqnarray*}\Omega_\mathcal{S}(X,Y)&=&-i\sum_{\alpha=1}^{m-1}\varepsilon_\alpha\cdot \langle \Omega(X,Y)\bullet
e_\alpha, Je_\alpha\rangle\\
&=&-i\sum_{\alpha=1}^{m-1}\varepsilon_\alpha\cdot \langle \Omega_0(X,Y)\bullet
e_\alpha, Je_\alpha\rangle\\
&=&-i\sum_{\alpha=1}^{m-1}\varepsilon_\alpha\cdot \eW(X,Y,e_\alpha,J_{ss}e_\alpha)\\
&=&-i\sum_{\alpha=1}^{m-1}\varepsilon_\alpha\cdot \big(\ \eW(X,J_{ss}e_\alpha,e_\alpha,Y)\ -\ \eW(X,e_\alpha,J_{ss}e_\alpha,Y)\ \big)\\
&=&-i\sum_{\alpha=1}^{m-1}\varepsilon_\alpha\cdot \big(\ \langle\Omega_0(X,J_{ss}e_\alpha)Je_\alpha,JY\rangle\ +\ 
\langle\Omega_0(X,e_\alpha)e_\alpha,JY\rangle\ \big)
\end{eqnarray*}
for all $X,Y\in TM$.

Let us assume now that one of the vectors $X$ and $Y$ is an element of $W(J,g)$. 
In fact, we can assume that $Y\in W(J,g)$ and $X\in TM$ is arbitrary. We set $u_{2\alpha-1}:= e_\alpha$ and $u_{2\alpha}:= 
J_{ss}e_\alpha$ for $\alpha=1,\ldots,m-1$.
With the remark from the beginning of our proof we 
obtain
\begin{eqnarray*}\Omega_\mathcal{S}(X,Y)&=&
\ -i\sum_{\alpha=1}^{m-1}\varepsilon_\alpha\cdot 
\big(\ \eW(X,J_{ss}e_\alpha,J_{ss}e_\alpha,J_{ss}Y)\ +\ \eW(X,e_\alpha,e_\alpha,J_{ss}Y)\ \big)\\[4mm]
&=&\ -\ i\sum_{i=1}^{n-2}\varepsilon_\alpha\cdot\eW(X,u_i,u_i,J_{ss}Y)\\
&&\qquad\quad +\quad \frac{i}{1+J_c^2}\cdot\left(\ \eW(X,j,J_\eta^\flat,J_{ss}Y)\ +\ \eW(X,J_\eta^\flat,j,J_{ss}Y)\ \right)\\[4mm]
&=&\ -i\cdot tr_g^{23} \eW(X,\ \cdot\ ,\ \cdot\ ,J_{ss}Y)\\[4mm]
&=&\ -i\cdot tr_g^{14} \eW(\ \cdot\ ,X,J_{ss}Y,\ \cdot\ )\\[3mm] 
&=&\quad 0\ .
\end{eqnarray*}
This shows that $\iota_Y\Omega_\eS=0$ for all $Y\in W(J,g)$. 
Since the orthogonal complement of $W(J,g)$ in $TM$ has dimension $2$ and spans 
together with $W(J,g)$ the tangent space $TM$, 
we can conclude 
that the only 
possible non-vanishing component of the curvature
on $\mathcal{S}$ is $\Omega_\eS(j,J_\eta^\flat)$. 
However, we know already that $\iota_j\Omega_\eS=0$, i.e., the latter component of $\Omega_\eS$
vanishes as well.	   
\hfill$\Box$\\

The main input for the proof of Lemma \ref{LM2} is the normalisation of 
$\omega_{nor}$ by the condition $\partial^*\circ \kappa=0$ on the conformal curvature (which  
is achieved for any space $(M,c)$). 
The proof uses then at a crucial point the generalised Bianchi identity, for which we introduce an
arbitrary metric in order to decompose standard tractors into densities and tangent vectors. 
Without the introduction of a metric it seems
not possible to apply the Bianchi identity to the tractor curvature.    
We note further that Lemma \ref{LM2} shows that the local complex tractor volume form $\varrho$ is 
parallel with respect to $\nabla^{nor}$, i.e.,
$\nabla^{nor}\varrho=0$. This property also implies the local existence of 
a conformal Killing  spinor on $(M,c)$ (cf. \cite{Bau99}). Lemma \ref{LM2} is the main ingredient 
for the proof of our reduction conjecture on the conformal holonomy.

\begin{THEO} \label{TH1} Let $(M,c)$ be a space of dimension $n=2m\geq 4$ with conformal structure $c$ of signature $(2p-1,2q-1)$
such that 
the conformal holonomy group $Hol(\mathcal{T})$ is contained in the unitary group $\U(p,q)$.
Then
\begin{enumerate}
\item\label{z2}
the holonomy algebra $\frak{hol}(\eT)$ of the canonical connection $\omega_{nor}$ is a subalgebra 
of the special unitary algebra $\frak{su}(p,q)$.
\item\label{z3}
If, in addition, the space $M$ is simply connected then the holonomy group $Hol(\mathcal{T})$ of $\omega_{nor}$ is contained in 
the special unitary group $\SU(p,q)$. 
\end{enumerate}
\end{THEO}

{\bf Proof.} (\ref{z2}) The assumption of Theorem \ref{TH1} about the holonomy group $Hol(\mathcal{T})$ implies the existence of a 
$\nabla^{nor}$-parallel complex structure $J\in\Gamma(\mathcal{A}(M))$ on the standard tractor bundle $\mathcal{T}(M)$
over the space $(M,c)$. 
The statement of Lemma \ref{LM2} shows that the values $\Omega(X,Y)\in\mathcal{A}(M)$ of the 
tractor curvature have vanishing complex trace for all $X,Y\in TM$. It follows that the curvature form $\Omega$ on 
the principal fibre bundle $\tilde{\mathcal{G}}(M)$ takes values only in $\frak{su}(p,q)\subset \frak{g}$.
Since the special unitary algebra $\frak{su}(p,q)$ is an ideal in $\frak{u}(p,q)$, the Ambrose-Singer Theorem proves 
that the holonomy algebra $\frak{hol}(\eT)$ of the canonical connection $\omega_{nor}$ on $\tilde{\mathcal{G}}(M)$
is contained in $\frak{su}(p,q)$. 

(\ref{z3}) In general, the holonomy group $Hol(\mathcal{S},\nabla^{nor})$ of 
the canonical complex line tractor bundle is a closed Lie subgroup
of $\U(1)$. Here, since $\mathcal{S}$ is locally flat by 
Lemma \ref{LM2}, $Hol(\mathcal{S},\nabla^{nor})$ is a discrete subgroup of $\U(1)$.
With the assumption that $M$ is simply connected it follows that 
$Hol(\mathcal{S},\nabla^{nor})$ is trivial, i.e., the complex line bundle 
$\mathcal{S}$ is globally flat and gives rise to a parallel complex tractor volume form on $M$. 
This proves that the conformal holonomy $Hol(\mathcal{T})$
is contained in $\SU(p,q)$. 
\hfill$\Box$\\

\section{The argument using Sparling's characterisation} \label{ab6} 

An integrable  CR-structure of hypersurface type on a manifold $N$ of odd dimension $2m-1\geq 3$ 
with signature $(p-1,q-1)$, $p+q=m+1$, of the Levi form gives rise to a $|2|$-graded 
parabolic geometry with structure group $\PSU(p,q)$. 
The classical Fefferman construction assigns to any such CR-structure a $S^1$-principal fibre 
bundle $M$ over $N$ of total dimension $2m$,
which is in a natural way equipped with a conformal structure $c$ of signature $(2p-1,2q-1)$. 
To be more precise, the explicit construction of $c$ is usually achieved by choosing a pseudo-Hermitian structure 
$\theta$ on $N$, which then gives rise to a particular metric $g$
on $M$. This metric is called the  Fefferman metric (corresponding to $\theta$) and 
it turns out that the conformal class of $g$ does not depend on the choice of pseudo-Hermitian form, 
but only on the CR-structure of $N$.
We do not aim either to introduce CR-geometry nor do we explain the Fefferman construction, in detail. 
For explanations on these subjects 
we refer e.g. to \cite{Fef76}, \cite{Lee86} or \cite{Cap02}. Instead, what is important for us is a characterisation of 
Fefferman metrics with help of a certain Killing vector field.
This is the so-called Sparling's characterisation of Fefferman metrics (cf. \cite{Spa85}, \cite{Gra87}).

\begin{THEO} \label{TH2} (Sparling's characterisation) (cf. \cite{Spa85}, \cite{Gra87}) Let $(M^n,g)$ 
be a pseudo-Riemannian space of dimension $n\geq 4$ and signature $(2p-1,2q-1)$.
Suppose that $g$ admits a Killing vector $j$ (i.e., $\mathcal{L}_jg=0$) such that
\begin{enumerate}
\item\label{z4}
$g(j,j)=0$, i.e., $j$ is lightlike,
\item
$\iota_j\eW=0$ and $\iota_j\mathcal{C}=0$, 
\item\label{z5}
$Ric(j,j) > 0$ on $M$.
\end{enumerate}
Then the metric $g$ is locally isometric to the Fefferman metric of some integrable CR-space $N$ of 
hypersurface type with signature $(p,q)$
and  dimension $n-1$.

On the other side, any Fefferman metric of an integrable CR-space $N$ of hypersurface type admits a Killing vector field $j$ 
satisfying conditions (\ref{z4}) to (\ref{z5}).
\end{THEO}

The Sparling's characterisation Theorem is suitable to prove Theorem \ref{TH1}.\\

{\bf Alternative Proof of Theorem \ref{TH1}.} Under the assumptions of 
Theorem \ref{TH1} there exists a $\nabla^{nor}$-parallel complex structure 
$J\in\Gamma(\mathcal{A}(M))$.
From Lemma \ref{LM0}, it follows that $j:=\Pi(J)$ is a conformal Killing vector on $(M,c)$. We choose an arbitrary metric $g$ in $c$.
It holds  $\mathcal{L}_jg=\frac{2}{n}div^gj\cdot g$.
Since $J$ is a complex structure, it follows from Lemma \ref{LM1} that the vector 
field $j$ is lightlike and admits no zeros on $M$. This implies that 
the partial differential equation
\[ j(\phi)\ =\ \frac{1}{n}\cdot div^gj\]
admits locally always a solution for $\phi$. 

We choose some local solution $\phi$ and proceed for the moment with local considerations.
We set $f:=e^{2\phi}g$ and  then 
the vector field $j$ is a Killing vector with respect to the metric $f$, i.e., $\mathcal{L}_jf=0$. 
From Lemma \ref{LM0}, it follows also that $j$ annihilates the conformal curvature $\Omega$. With respect to the metric $f$ this means
$\iota_j\eW=0$ and $\iota_j\mathcal{C}=0$. Moreover, it holds
\[
div^fj=0\qquad\mbox{and}\qquad\mathcal{D}^f(j)\ =\ \frac{1}{n-2}\left(\ -Ric^f(j,\cdot)\ 
+\ \frac{scal^f}{2(n-1)}\cdot f(j,\cdot)\ \right) \ . 
\]  
Since $f(j,\mathcal{D}^f(j))=-1-(div^fj)^2$, these identities  prove that \[Ric^f(j,j)\ =\ n-2\ >\ 0\ .\]
The results so far show that we can apply Theorem \ref{TH2}. We conclude that $f$ is (locally) isometric to the Fefferman metric 
of some integrable CR-space of hypersurface type.
  
Similar as in conformal geometry it is the case that a CR-structure on some space $N$ 
is uniquely determined by a normal Cartan connection on some principal 
fibre bundle. The curvature of this unique normal Cartan connection takes values 
in $\frak{su}(p,q)$, which is the Lie algebra of the structure group of CR-geometry. 
Moreover, it is known in the classical Fefferman construction for integrable CR-spaces 
that the lift along the $S^1$-fibration of the unique normal Cartan connection
of CR-geometry gives rise to the unique normal  Cartan connection of
the Fefferman conformal class. In our situation this argument says that the conformal curvature $\Omega$ of $c=[f]$
takes values only in $\frak{su}(p,q)$.

The latter argument is rather vague, since we use some construction from CR-tractor calculus, 
which we do not develop here (cf. \cite{Cap02}).
However, we can replace this argument by the following one. It is known that any Fefferman metric admits at least locally
a certain conformal Killing spinor (cf. \cite{Bau99}). This fact implies again that the conformal curvature $\Omega$ of $c=[f]$
has values solely in $\frak{su}(p,q)$.

The latter statement about the curvature is locally true everywhere on $(M,c)$ and we conclude that the conformal curvature $\Omega$
has no complex trace on $(M,c)$. Again, as in the proof of Theorem \ref{TH1}, 
we can apply the Ambrose-Singer Theorem to reason that the 
conformal holonomy algebra $\frak{hol}(\mathcal{T})$ of $(M,c)$ is contained in $\frak{su}(p,q)$. 
This proves the first statement of Theorem \ref{TH1}, which implies as before with 
further assumption the second statement
on the holonomy group $Hol(\mathcal{T})$. \hfill $\Box$\\ 

Sparling's characterisation is not only useful to prove the result of Theorem \ref{TH1}. In addition, it  
gives rise to a geometric characterisation and construction principle for spaces whose conformal holonomy sits in $\U(p,q)$.  

\begin{COR} \label{COR100} Let $(M,c)$ be a space with conformal structure of signature $(2p-1,2q-1)$ and conformal holonomy
group contained in $\U(p,q)$. Then the conformal class $c$ is locally the Fefferman conformal class of
some integrable CR-space of hypersurface type of signature $(p,q)$.  
\end{COR}


\end{sloppypar}

\begin{thebibliography}{1111111}

\bibitem[Arm05]{Arm05} S. Armstrong. {\it Conformal Holonomy: A Classification}.
e-print: arXiv:math.DG/ 0503388, (2005).

\bibitem[BEG94]{BEG94} T.N. Bailey. M.G. Eastwood. A.R. Gover. 
{\it Thomas's structure bundle for conformal, projective and related structures}.
Rocky Mountain J. Math.  24  (1994),  no. 4, 1191--1217.

\bibitem[Bau99]{Bau99} H. Baum. {\it Lorentzian twistor spinors and CR-geometry}.
Differential Geom. Appl. 11 (1999), no. 1, 69--96.


\bibitem[Cap02]{Cap02} A. Cap. {\it
Parabolic geometries, CR-tractors, and the Fefferman construction}. Differential Geom. Appl. 17 (2002), 123--138.


\bibitem[Cap05]{Cap05} A. Cap. {\it Infinitesimal Automorphisms and Deformations of Parabolic Geometries}.
ESI e-Preprint 1684, Vienna, 2005.

\bibitem[CG03]{CG03} A. Cap. {\it Standard tractors and the conformal ambient metric construction}. 
Ann. Global Anal. Geom. 24, 3 (2003) 231-259.



\bibitem[CSS97a]{CSS97a} A. Cap. J. Slovak. V. Soucek. {\it
Invariant Operators on Manifolds with Almost Hermitian Symmetric Structures, I.
Invariant Differentiation}. Acta Math. Univ. Comenian.
66 no. 1 (1997) 33-69.


\bibitem[CSS97b]{CSS97b} A. Cap. J. Slovak. V. Soucek. {\it
Invariant Operators on Manifolds with Almost Hermitian Symmetric Structures, II. Normal Cartan Connections}.
Acta Math. Univ. Comenian. 66 no. 2 (1997), 203-220.

\bibitem[CSS01]{CSS01}
A. Cap. J. Slovak. V. Soucek. {\it Bernstein-Gelfand-Gelfand sequences}.  Ann. of Math. (2)  154
(2001),  no. 1, 97--113.


\bibitem[Fef76]{Fef76} Ch. Fefferman. {\it Monge-Ampere equations,
the Bergman kernel, and geometry of pseudoconvex domains}. Ann. Math. 103
(1976), 395-416.


\bibitem[Gra87]{Gra87} C.R. Graham. {\it On Sparling's characterization of Fefferman metrics}.
Amer. J. Math. 109 (1987), no. 5, 853--874.

\bibitem[Kob72]{Kob72}
Sh. Kobayashi. {\it Transformation groups in differential geometry}. Ergebnisse der Mathematik und ihrer
Grenzgebiete,
Band 70. Springer-Verlag, New York-Heidelberg, 1972.



\bibitem[Lee86]{Lee86} J. M. Lee. {\it The Fefferman metric and
pseudo-Hermitian invariants}.  Trans. Amer. Math. Soc.  296  (1986),  no.
1, 411--429.


\bibitem[Lei04]{Lei04} F. Leitner. {\it Conformal holonomy of bi-invariant metrics}.
e-print: arXiv:math.DG/ 0406299, (2004).



\bibitem[Lei05a]{Lei05a} F. Leitner. {\it Conformal Killing forms
with normalisation condition}. Rend. Circ. Mat. Palermo Suppl. ser. II. (2005), no. 75, 279-292.

\bibitem[Lei05b]{Lei05b} F. Leitner. {\it About complex structures in conformal tractor calculus}.
e-print: arXiv: math.DG/0510637, (2005). 


\bibitem[Spa85]{Spa85} G.A.J. Sparling. {\it Twistor theory and the 
characterisation of Fefferman's conformal structures.} Preprint Univ.
Pittsburg, 1985.




\end{thebibliography}
\end{document}